\documentclass{article}
\usepackage{graphics}
\usepackage{latexsym}
\usepackage{amssymb}
\usepackage{amsmath}
\usepackage{amsthm}
\usepackage{epsfig}
\renewcommand{\deg}{^{\circ}}

\newtheorem{definition}{Definition}
\usepackage{epstopdf}

\begin{document}
\begin{center}
\vspace*{-0.8in}
\begin{minipage}{\textwidth}
{\bf\Huge \textsf{Forcing Nonperiodicity \\ \vspace{1pt} \\ With a Single Tile}} 

\vspace{12pt}
{\bf\large \textsf{Joshua E.~S.~Socolar$^1$ and Joan~M.~Taylor$^2$}} \\
\textsf{$^1$Physics Department and Center for Nonlinear and Complex Systems
 Duke University, Durham, NC 27708} \\
 \textsf{$^2$P.O. Box U91, Burnie, Tas. 7320 Australia}
\end{minipage}
\end{center}

\vspace{18pt}
It is easy to create nonperiodic tesselations of the plane composed of one or a few types of tiles.  In most cases, however, the tiles employed can also be used to create simpler, periodic patterns.  It is much more difficult to find shapes, or ``prototiles,'' that can fill space {\em only} by making a nonperiodic structure.   We say that such sets are {\it aperiodic}, or that they ``force'' nonperiodicity, and there are many open questions about what types of structure can be forced and the prototiles required.  In this article we discuss recent progress on the fundamental problem of forcing nonperiodicity using a {\it single} prototile, jokingly called an {\it einstein}  (a German pun on ``one stone'').   A new example we found~\cite{SocolarTaylorJCT}  shows one way in which an einstein can work and highlights several issues that arise in posing the problem precisely.

One motivating factor in the search for an einstein comes from condensed matter physics.  Local rules for how tiles fit together may represent the energetics of a physical system, which could support self--assembly into an ordered but nonperiodic structure.  The discovery of icosahedral and decagonal phases of metallic alloys, in which the atomic structure shares the essential structure of the Penrose tilings, has opened our eyes to the fact that nonperiodic materials can indeed form spontaneously~\cite{SteinhardtBook,janot}.  In materials physics applications, where the tiles may represent clusters of many atoms or larger building blocks, the tiles can have complex shapes or markings that determine how they may be joined.  Finding a single shape that can do the job may make the physical realization of such a material easier.

The first example showing that it is possible to force nonperiodicity was Berger's set of $20,426$ distinct prototiles \cite{Berger}.   Aperiodic sets with just two prototiles were subsequently discovered, the most famous being the Penrose tiles~\cite{Penrose74}, nicely described by Martin Gardner~\cite{Gardner77}.   Candidates with einstein--like features have been presented before, but there is no precise definition of the einstein problem, and several candidates that could be argued to qualify have not passed the consensus ``I know it when I see it'' test.  There are several issues involved, including the specification of what counts as nonperiodic, what characteristics make for a valid prototile, and what form the local rules must take. 

\begin{figure}[htbp]
\begin{center}
\includegraphics[scale=0.19]{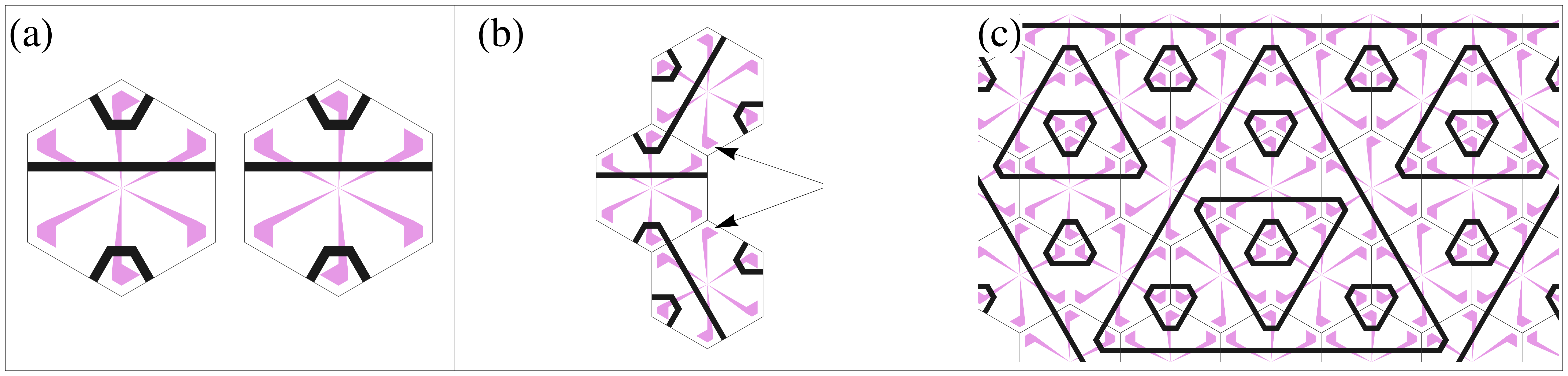}\hspace{36pt}
\caption{The hexagonal prototile and its mirror image with color matching rules.  (a) The two tiles  are related by reflection about a vertical line.  (b)  Adjacent tiles must form continuous black stripes.  Flag decorations at opposite ends of a tile edge, such as the indicated flags at opposite ends of the vertical edge,  must point in the same direction.   (c) A portion of an infinite tiling that respects the matching rules.}
\label{fig2Dmirrors}
\end{center}
\end{figure}

We recently showed that the prototile in Figure~\ref{fig2Dmirrors} is an einstein and determined a number of remarkable properties of the tilings it forces.  Ref.~\cite{SocolarTaylorJCT} contains two proofs of the forced nonperiodicity along with derivations of several intriguing properties of the tiling (including a surprising connection to the regular paperfolding sequence~\cite{paperfoldingMI}).  In working out the properties of the forced limit--periodic structure and searching for different ways of encoding the information about how the tiles must fit together, we were led to a series of questions about how the einstein problem should be posed.   In the present paper, we discuss the definitions of the terms {\it local matching rules} and  {\it tile}, and we propose a new definition of {\it nonperiodic} that emphasizes distinctions that to our knowledge have not been made explicit before.  We use our hexagonal prototile throughout to clarify key points, including a 3D version for which the shape alone is sufficient to force nonperiodicity.   We also describe some of the more intriguing aspects of the forced structure, some of which require further study.  Articulating the criteria satisfied by our hexagonal prototile reveals the sense in which it is the ``best'' einstein currently known and delineates a precise problem that remains open.

The newly discovered tile, shown together with its mirror image in Figure~\ref{fig2Dmirrors}(a), is a regular hexagon with markings that determine how neighboring tiles must be oriented.  Adjacent tiles must form continuous black stripes and flag decorations at opposite ends of each tile edge must point in the same direction.   (The arrows in (b) point to the two flags at opposite end of a vertical tile edge.)  Each tile in (c) is a rotation and/or reflection of the single prototile and the only way to fill space while obeying the rules everywhere is to form a nonperiodic, hierarchical extension of the pattern in (c). 

\section*{Defining the einstein}

Two constructions that could conceivably be counted as einsteins were discovered in 1995.  A single prototile that forces a pattern of the Penrose type was presented by Gummelt  (with a complementary proof by Steinhardt and Jeong)~\cite{Gummelt96,Jeong96}.  But in this case tiles are allowed to overlap and the covering of the space is not uniform.  For this reason the prototile is not considered to be an einstein.

\begin{figure}[htbp]
\begin{center}
\includegraphics[scale=0.19]{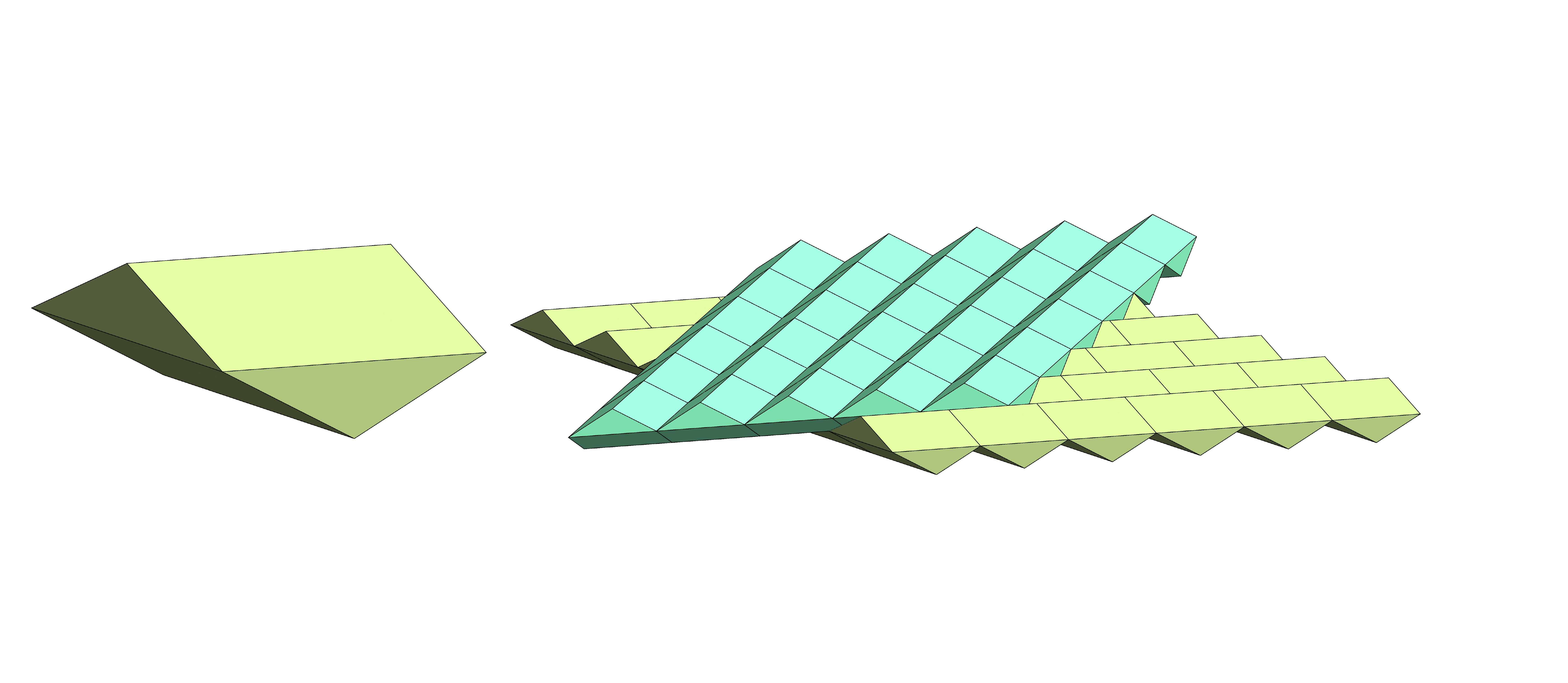}\hspace{36pt}
\caption{The SCD prototile and the space--filling tiling it forces.}
\label{fig:scd}
\end{center}
\end{figure}

The uniformly space--filling, three dimensional prototile of Figure~\ref{fig:scd}, a rhombic biprism, was exhibited by Schmitt, Conway and Danzer~\cite{SCD}.  
To fill space, one is forced to construct 2D periodic layers of tiles sharing triangular faces, with ridges running in the direction of one pair of rhombus edges on top and the other pair below.  The layers are then stacked such that each is rotated by an angle $\phi$ with respect to the one below it, where $\phi$ is the acute angle of the rhombic base.  Any choice of $\phi$ other than integer multiples of $\pi/3$ or $\pi/4$ produces a tiling that is not periodic, and certain choices permit a tiling in which the number of nearest neighbor environments is finite, so that the prototile can be endowed with bumps and nicks in a way that locks the relative positions of adjacent layers.

Again, however, the universal reaction was ``This is not really what we are looking for.''   The nonperiodicity of the tiling does not seem mysterious enough to count; one can immediately grasp the global structure of simple 2D periodic lattices stacked with a twist in the third dimension.  We seek structures with long--range correlations that are not immediately evident from the examination of a single tile.  Goodman--Strauss emphasizes this point and suggests a classification that distinguishes the SCD structure from the Robinson tiling and Penrose tiling~\cite{GoodmanStraussQs}.  Goodman--Strauss calls the SCD tile {\it weakly aperiodic} because it admits a tiling with a cyclic group of symmetries involving finite (and nonzero) translations, in this case the screw operations along the twist axis.  We emphasize here an additional feature of the SCD tiling that also weakens the sense in which it may be called nonperiodic:  every individual tile in the tiling is a unit cell of a periodic 2D layer.  In fact, for the cases with a finite number of nearest neighbor environments, any finite stack of layers is periodic in the two transverse directions, a point that may be important for physical applications.  The presence of infinite periodic substructures within the pattern suggests that a formal definition be developed to distinguish different degrees of periodicity or nonperiodicity.  

\subsection*{A definition of ``nonperiodic''}
We offer here a new classification scheme, based on the notion of a ``partial translational symmetry,'' that we believe captures the shared intuitive notion of a nontrivial, nonrepeating pattern.  
A partial translational symmetry is an operation that maps some subset of a full pattern into itself.  
\begin{definition}  Let ${\cal T}$ be a infinite set of tiles in ${\cal R}^N$.  A {\em partial translational symmetry} of ${\cal T}$ is an operation of the form $\{\vec{x} \rightarrow \vec{x} + {\bf R}\cdot\vec{x} +  \vec{e}\}$  that acts on some infinite subset of ${\cal T}$ and leaves it invariant,  where ${\bf R}$ is a rotation matrix and $\vec{e}\in{\cal R}^N$ is a constant, nonzero, displacement vector.  The magnitude of $\vec{e}$ is the {\em spacing} of the partial translational symmetry.
\end{definition}
We say that a tile {\em participates} in a partial translational symmetry if it is a member of the subset of tiles that is left invariant by the symmetry operation.  The tiles participating in a given partial translational symmetry need not form a connected region.

The SCD tiling has many partial translational symmetries.  The tiles in any one layer form a subset that is invariant under a 2D lattice of translations (for which ${\bf R}$ is the identity).  In addition, every tile is an element of a subset that is invariant under a screw operation (${\bf R}$ being a rotation about the stacking axis) that maps one layer into the next.  

\begin{figure}
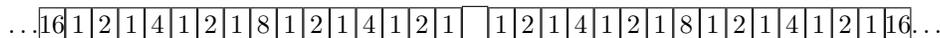

$\ldots$\framebox[10pt]{16}\framebox[10pt]{1}\framebox[10pt]{2}\framebox[10pt]{1}\framebox[10pt]{4}\framebox[10pt]{1}\framebox[10pt]{2}\framebox[10pt]{1}\framebox[10pt]{8}\framebox[10pt]{1}\framebox[10pt]{2}\framebox[10pt]{1}\framebox[10pt]{4}\framebox[10pt]{1}\framebox[10pt]{2}\framebox[10pt]{1}\framebox[10pt]{\rule{0pt}{7pt}\ }\framebox[10pt]{1}\framebox[10pt]{2}\framebox[10pt]{1}\framebox[10pt]{4}\framebox[10pt]{1}\framebox[10pt]{2}\framebox[10pt]{1}\framebox[10pt]{8}\framebox[10pt]{1}\framebox[10pt]{2}\framebox[10pt]{1}\framebox[10pt]{4}\framebox[10pt]{1}\framebox[10pt]{2}\framebox[10pt]{1}\framebox[10pt]{16}$\ldots$ 
\caption{A limit--periodic pattern with partial translation symmetries.}
\label{fig:gcd}
\end{figure}

Another simple example of a system with nontrivial partial translational symmetries will help clarify the sense in which our new tiling is nonperiodic.   Figure~\ref{fig:gcd} shows a 1D tiling in which tiles have unit width and the number on a tile indicates its type.
The tiling extends infinitely in both directions and the spacing between the nearest tiles of the same type is just twice their numerical value; type $x$ occurs periodically at positions $x(1+2 n)$ for integer $n$.   (The blank element at the center is not repeated at any finite distance.)  Thus, for each value $y$ occurring in the sequence, the subset of tiles with values less than or equal to $y$ is invariant under translation by $2 y$.  For example, under translation by $8$, the set of all tiles with values less than or equal to $4$ is invariant, though the set of remaining tiles is not.

\begin{definition}  The elements of a set of partial translational symmetries are {\em independent} if and only if it is not possible to express the displacement $\vec{e}$ for any one of them as an integer linear combination of the others.  
\end{definition}
In an ordinary $N$--dimensional periodic tiling, every tile participates in $N$ independent partial translational symmetries, the displacements $\vec{e}_i$, with $i = 1,\ldots N$, being the basis vectors of the lattice of translations that leave the tiling invariant. 

\begin{definition}  Let the number of independent partial translational symmetries with spacing less than $r$ that a given tile $t$ participates in be denoted $S_t(r)$.  A tiling in $\Re^N$ is {\em nonperiodic} if for any finite $r$, the fraction of tiles with $S_t(r)$ strictly less than $N$ is finite.
\end{definition}
By this definition, the decagonal quasicrystal structures~\cite{Kortan90}, which are periodic stackings of quasicrystalline layers, are nonperiodic because the only independent partial translational symmetry containing any given tile is the one corresponding to the periodic stacking direction.  (There could conceivably be additional screw operations if the tiling has an axis of complete 5-fold rotational symmetry, but they all have the same $\vec{e}$ or integer multiples of it.)

But the SCD tiling is {\em not} nonperiodic by this definition.  As noted above every tile is a member of a subset that is invariant under two distinct translation operations in the plane and a subset that is invariant under a screw operation along the stacking direction.  To emphasize the high degree of partial translational symmetry in the SCD structure, we might classify it as {\em heterogeneously periodic}; ``heterogeneously'' because, unlike familiar periodic structures, the $\vec{e}$'s associated with the symmetries of tiles in different layers are not all the same.   Heterogeneously periodic is not synonymous with weakly nonperiodic; some weakly nonperiodic tilings, such as the stacked layers that form decagonal quasicrystals, are nonperiodic by our definition (and by common usage in the physics community).
  
Many tilings satisfy a more stringent criterion:
\begin{definition}  A tiling is {\em maximally nonperiodic} if and only if it contains no partial translational symmetries.
\end{definition}
The Penrose tilings are maximally nonperiodic, as are many tilings generated by substitution rules or projections from higher dimensional periodic lattices onto incommensurate subspaces.  

Following Gr\"{u}nbaum and Shephard, we adopt the following terminology.  
\begin{definition}  A prototile is (maximally) {\em aperiodic} if the only space--filling tilings that can be composed from it are (maximally) nonperiodic.
\end{definition}
\noindent That is, we call a single prototile ``aperiodic'' if it can be used to tile an entire space with no overlaps but  only in a pattern that is nonperiodic by the above definition.  

\subsection*{The question of allowable matching rules}
Whether or not an aperiodic prototile exists could hinge on whether one requires that the tile be a simply connected domain, whether tiles are allowed to overlap so that some parts of space are doubly covered, whether the rules must be encoded by tile shape alone as opposed to color--matching rules, and whether mirror--image sets are considered to count as a single prototile or not.   A highly restrictive definition of an einstein would demand the following properties:
\begin{description}
\item[Rotations only:]  Reflections of a chiral tile are {\em not} allowed;
\item[Simply connected tiles:]  The prototile is a simply connected domain (a topological disk in 2D);
\item[Shape alone:]  All configurations of tiles that do not contain overlaps are permitted without regard to any colored markings;
\end{description}

The first condition makes what some may view as an arbitrary distinction between rotations and reflections.  Nevertheless, we see a meaningful distinction between  cases where the tiles could all be manufactured from a single physical mold and cases where a second, mirror image mold must be built.

The restriction to simply connected prototiles is consistent with the intuitive notion of a tile as a thin, rigid piece of material, as is used in mosaics or floor tiling.  On the other hand, there is no obvious reason to insist that a tile cannot be composed of set of disconnected domains with fixed relative positions~\cite{GS}, and in fact certain types of color matching rules that cannot be enforced by the shape alone using a simply connected 2D prototile can be enforced using a tile consisting of disconnected pieces~\cite{Socolar07}.  Gr\"{u}nbaum and Shephard make a further distinction between tiles with cutpoints (where regions are connected only through a vertex) and tiles with entirely disconnected regions~\cite{GS}.  From a materials physics perspective, tiles may represent complex atomic configurations with low energy and these may conceivably interpenetrate in ways that could not be represented by simply connected tiles.

The ``shape alone'' condition requires further comment, as there are several kinds of rules that cannot be encoded in the geometry of the prototile.  
\begin{description}
\item[Colors required:]  
Instead of bumps and nicks, the rules that force relative orientations of nearby tiles can be encoded as a colored decoration of the prototiles together with rules about how colors must match.  
Not every rule enforced by color matching can be implemented through shape alone without increasing the number of prototiles.  A classic example is Ammann's aperiodic set $A5$ (a square and a $45\deg$ rhombus), where rules for how tiles must join at vertices may be implemented either through constraints on colored decorations around the vertex or by introduction of a new tile that must fit at each vertex~\cite{GS}.  Another example is the hexagonal parquet tile of \cite{Socolar07}, for which the color rules for tile edges (either red or black can match black, but red cannot match red) could be implemented by introducing two new tiles that fit into notched edges.  
\item[Non--adjacent, but pairwise:]  A rule may specify the relative orientations of two tiles separated by some bounded distance but not sharing an edge.  In such a case, it is still possible to check whether the tiling satisfies the rules by examining only two tiles at a time, or, as physicists would say, by considering only pairwise interactions between tiles.
\item[Configuration atlases:]  The set of allowable configurations may be expressed as an atlas of allowed configurations within some ball of finite radius, but not be expressible as a set of pairwise constraints.  Examples include the trivial cases presented by Goodman-Strauss in which rectangular tiles are required to form pixellated versions of Robinson square tiles~\cite{GoodmanStraussQs} as well as the recent construction of Fletcher in which face-matching rules for an aperiodic set of 21 cubic prototiles are expressed as an atlas of allowable configurations for a single cubic prototile in which the 21 different tile types are encoded as 21 different orientations of a single tile~\cite{Fletcher}.
\end{description}

\section*{The einstein}
The discovery of the prototile and rules of Figure~\ref{fig2Dmirrors} was initiated by Taylor's observation that a single colored hexagon together with its mirror image could force a structure similar to the one forced by a set of 12 tiles (discovered later by Socolar and Goodman-Strauss) appearing on Socolar's web page.  She had been searching since 1993 for a superposition of matching rules to force nonperiodicity on the hexagon with black stripes, viewing it as an elementary version of the Penrose rhombi, which form a quasiperiodic hierarchy of overlapping, irregular hexagons.  Taylor's constructions were based heavily on a complex scheme for generating the tilings through a substitution procedure in which each tile is divided into smaller tiles that respect the same local rules.~\cite{substitutiontilingweb}  A note from Taylor requesting feedback led Socolar to refine the set of necessary local rules and construct a simple proof of aperiodicity, which initiated an extended collaboration conducted entirely by email between Tasmania and North Carolina.    The conceptual breakthroughs needed to resolve various subtle issues came about through repeated exchange of figures and discussion of details specific to these tilings.  There was no clearly generalizable strategy involved, though we hope that our results will lead by example to further discoveries.

The prototile and its mirror image are shown in Figure~\ref{fig2Dmirrors} as regular hexagons decorated with colors that encode rules constraining the relative orientation of nearby tiles.   There are two such constraints, or {\em matching rules}:  ({\bf R1})  the black stripes must be continuous across all edges in the tiling; and ({\bf R2}) the flags at the vertices of two tiles separated by a single tile edge must always point in the same direction.  The rules are illustrated in Figure~\ref{fig2Dmirrors}(b) and a portion of a tiling satisfying the rules is shown in Figure~\ref{fig2Dmirrors}(c).\footnote{We note that this tiling is similar in many respects to a tiling exhibited previously by Penrose~\cite{Penrose97}, but the two are not equivalent~\cite{SocolarTaylorJCT}. }

\begin{figure}[tb]
\begin{center}
\includegraphics[scale=0.19]{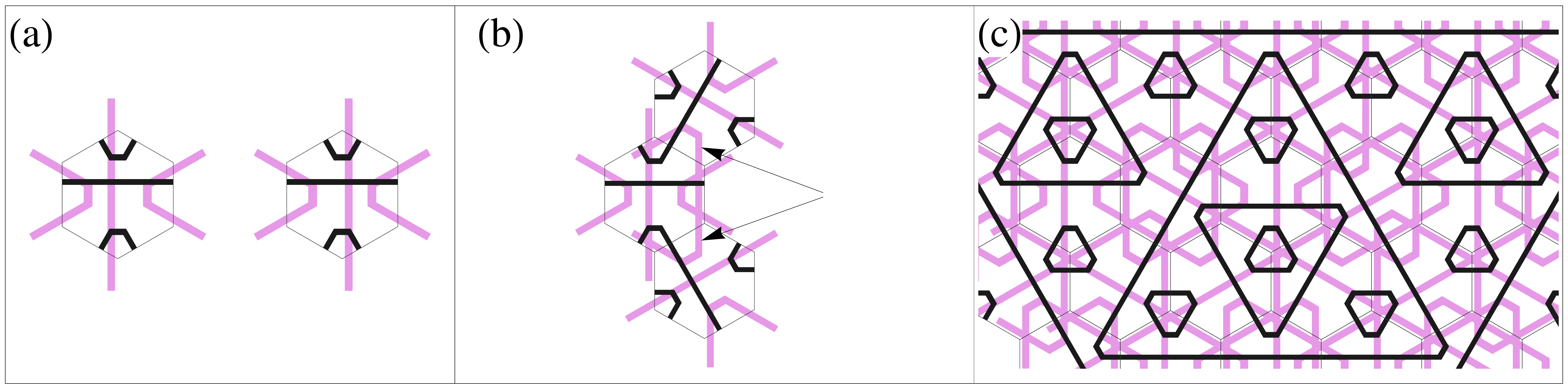}
\caption{Alternative coloring of the 2D tiles.  The arrows indicate stripes on next nearest neighbor tiles that must join to form a continuous line.}
\label{figblackpurple}
\end{center}
\end{figure}

An alternative way to represent the matching rules is to allow decorations that extend beyond the tile edges as shown in Figure~\ref{figblackpurple}.  {\bf R1} remains the same and {\bf R2} is now enforced by requiring that the purple stripes be continuous.  This makes it clear that {\bf R2} and {\bf R1} have identical geometric forms related by a scale factor of $\sqrt{3}$ and a rotation by $\pi/2$.  It is the relative positions of the long black and long purple stripe that distinguishes the different reflections of the prototile.  The purple stripes in Figure~\ref{figblackpurple} form hierarchical triangular structures just like the black stripes, but there are three interpenetrating purple structures.  (See Figure~\ref{figpartialtrans} below.)

The only space--filling tilings allowed by the 2D prototile of Figure~\ref{fig2Dmirrors} are nonperiodic.  The proof given in \cite{SocolarTaylorJCT} shows that the tiling forms an infinite hierarchy of interpenetrating  honeycomb lattices of black rings and the tiles in lattices with translational symmetry scales larger than $r$ cannot participate in partial translational symmetries with spacing smaller than $r$.  For any $r$, the density of tiles in larger scale lattices is clearly finite (nonzero), so the tiling is indeed nonperiodic by our definition.  Thus, assuming that we allow color matching rules and count mirror images as a single prototile, we have an einstein that requires only pairwise matching rules!

The partial translational symmetries of the black ring structure are immediately clear, but because of the interplay between that and the purple stripes, the partial translational symmetries of the full pattern are more difficult to locate.  Figure~\ref{figpartialtrans} displays a subset of them.  The shaded tiles are a motif that is repeated periodically to form a triangular lattice.  For visual clarity, the black stripe decoration is displayed for every tile and the purple stripe decoration of Figure~\ref{figblackpurple} is shown only for a subset of the tiles that form a triangular lattice consisting of one third of all of the tiles.  There is no partial translational symmetry in this tiling with a smaller spacing than the ones shown except for the special cases where partial translational symmetries occur along one infinite line in the tiling.  The proof involves analysis of the separate symmetries of the black stripe and purple stripe patterns.  We omit it here because it is not terribly illuminating and we have not yet solved the general problem of which tiles participate in partial translational symmetries with given spacings.

\begin{figure}[h]
\begin{center}
\includegraphics[scale=0.75]{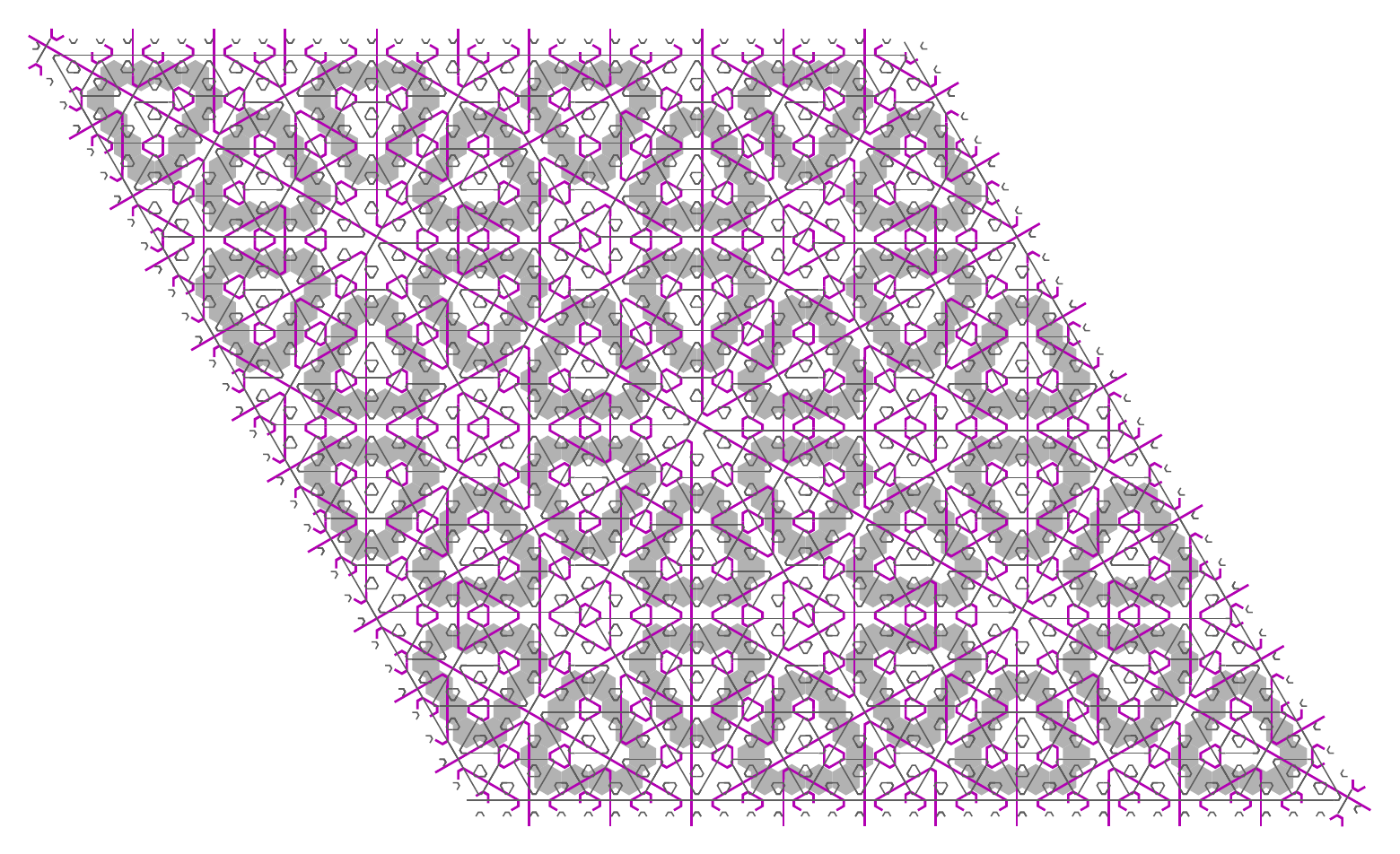}
\caption{The partial translational symmetry with the smallest spacing.  Clusters of 24 shaded tiles (two of each of the twelve tile orientations) are repeated throughout the tiling, forming a triangular lattice. Purple stripes are shown only for a subset of one third of the tiles.}
\label{figpartialtrans}
\end{center}
\end{figure}

\begin{figure}[tb]
\begin{center}
\includegraphics[scale=0.4]{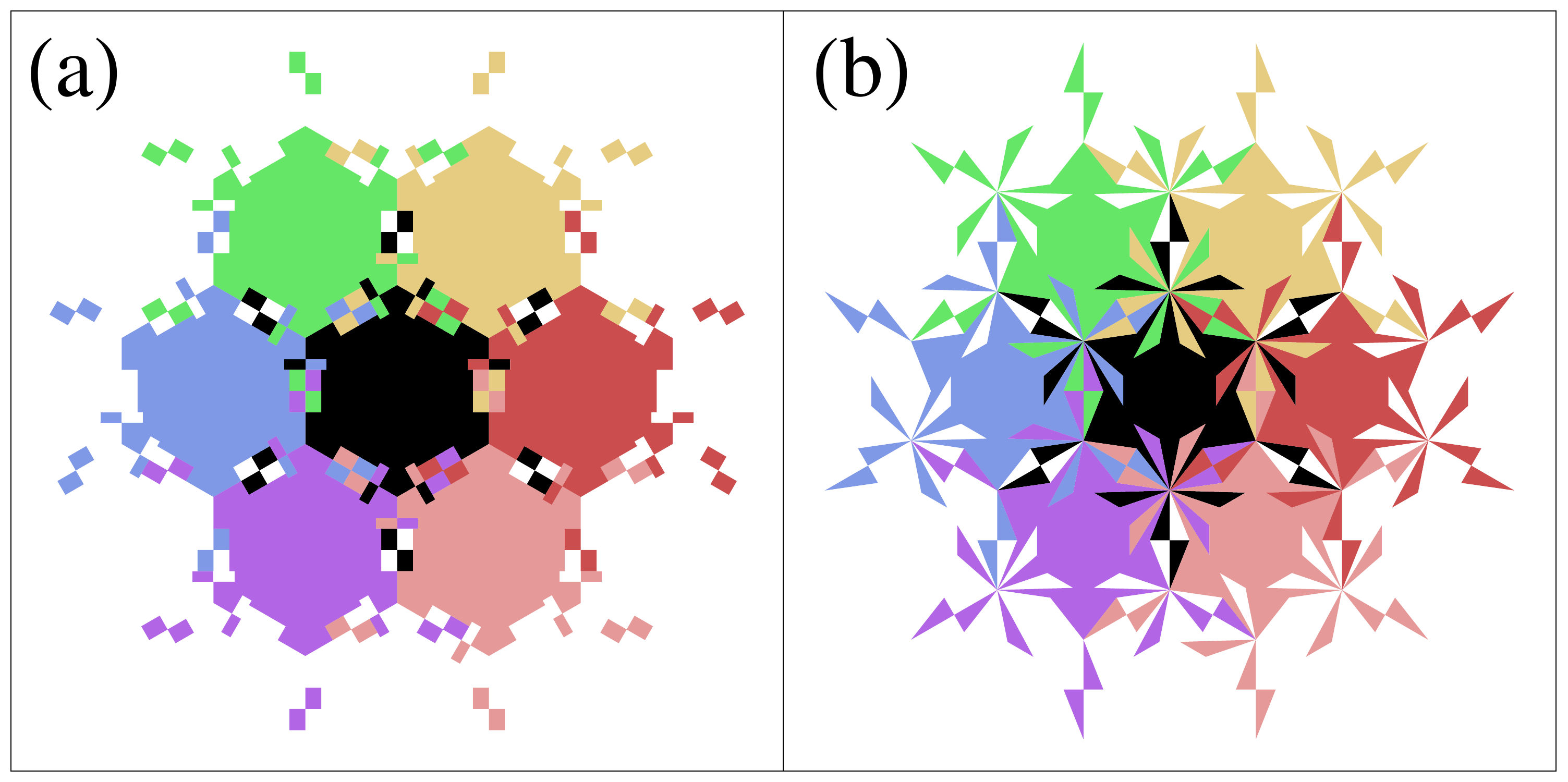}
\caption{(a) Enforcing the by shape alone with a disconnected 2D tile.  All the patches of a single color, taken together, form a single tile.  (b)  A deformation of the disconnected prototile in (a) to a prototile with cutpoints.}
\label{fig2Dmultcon}
\end{center}
\end{figure}

The matching rules {\bf R1} and {\bf R2} may appear to be unenforceable by shape alone.  {\bf R2} necessarily refers to tiles that are not in contact in the tiling and {\bf R1} cannot be implemented using only the shape of a single prototile and its mirror image.  Both of these obstacles can be overcome, however, if one relaxes the restriction that the prototile must be a simply connected shape.  Figure~\ref{fig2Dmultcon}(a) shows how the color--matching rules can be encoded in the shape of a single prototile that consists of several disconnected regions.  In the figure, all regions of the same color are considered to compose a single tile.  {\bf R1} is enforced by the small rectangles along the tile edges.  {\bf R2} is enforced by the pairs of larger rectangles located radially outward from each vertex.  The flag orientations are encoded in the chirality of these pairs.  Thus we have an einstein that does not require color matching rules!  Figure~\ref{fig2Dmultcon}(b) shows a deformation of the disconnected prototile to a prototile with cutpoints; i.e., a tile in which all the pieces are connected through vertices and tiles are allowed to overlap at those points.  For a beautiful rendering of this construction, see Araki's beetles~\cite{Araki}.

\begin{figure}[tb]
\begin{center}
\includegraphics[scale=0.40]{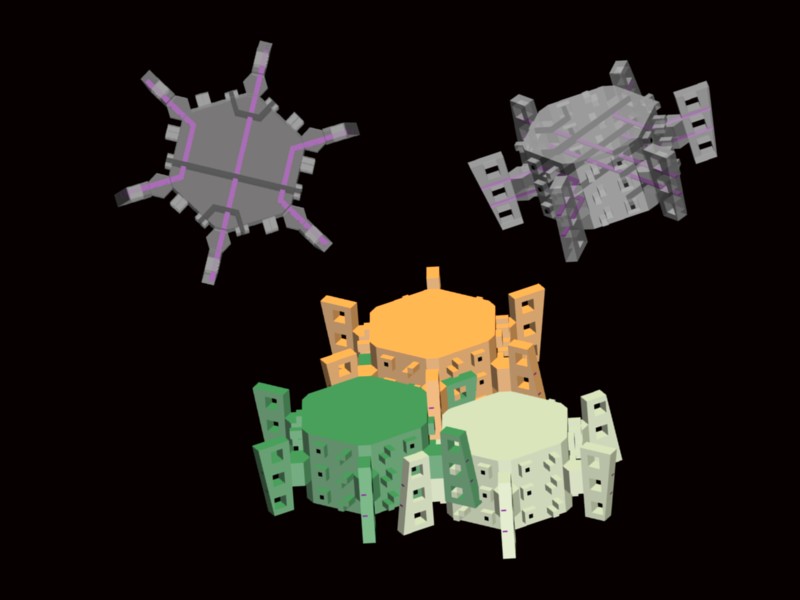}
\caption{A 3D einstein.  Colored bars are included only to clarify the relation of the 3D shape to the matching rules for the 2D tile.  A translucent tile is shown in two orientations to emphasize the relation of the colored bars to the shape. Three solid tiles that fit together properly are shown (with the ends of the colored bars visible).}
\label{fig3Dbrick}
\end{center}
\end{figure}

Whether you prefer to enforce the matching rules using colors or a disconnected prototile is a matter of taste.  Of course you may find both less than fully satisfying, in which case we can offer a third way out --- via escape to the third dimension.  The tiles of Figure~\ref{fig2Dmirrors} are related by reflection through a line in the 2D plane, but they can also be thought of as related by a rotation in 3D space of $180^{\circ}$ about that same line, suggesting that the two mirror image tiles be thought of as the front and back faces of a single 3D tile.  Such a tile is shown in Figure~\ref{fig3Dbrick}.

The colored bars running through the 3D tile are guides to the eye that display the black and purple stripe structure, but they are not required.  The continuity of the bars is enforced by the shape of the tile alone.  To see how, consider first the flag matching rule {\bf R2}.   To enforce this rule, we must have arms extending outward from the basic hexagonal prism to meet with the arms of next--nearest neighbor hexagons.  At each vertex of the hexagonal tiles, three arms must somehow pass through each other.  The tile shown in Figure~\ref{fig3Dbrick} solves this problem by allowing tiles to be staggered at three different heights.  The full tiling is divided into three triangular lattices of tiles, each of which contains tiles at one height.  The top faces of the tiles in the three different lattices are at heights $0$, $h/3$, and $2h/3$, where $h$ is the height (or thickness) of a tile.  The hexagonal blocks on each arm have thickness $h/3$, allowing the blocks from three crossing arms to make a full column.  The six arms on the prototile have outer faces that are tilted from the vertical in a pattern that encodes the chirality of the flags of the 2D tile.  Forming one triangular lattice requires that bevels of opposite type be joined, and hence that flags of opposite chirality match in accordance with {\bf R2}.

The small bumps on the tiles and the holes in the arms are arranged such that adjacent tiles can fit together if and only if the black stripes match up properly, as required by {\bf R1}.   The three square holes in each arm are positioned so that projections from the faces on neighboring tiles can meet with each other.  The holes are all the same; they do not themselves encode the positions of the black stripes.  Next, we create two types of plug that can be inserted into a hole.  One type consists of two square projections that fill opposite quadrants of the hole; the other type fills the entire hole but only to half its depth.  The two types are both invariant under rotation by $180^{\circ}$.  Two plugs of the same type can fit together to fill a hole, but plugs of different types cannot.  Finally, we place two columns of three plugs each on each of the large vertical faces of the main hexagonal portion of the tile.  Each column aligned with a black stripe has plugs of one type, and the other columns have plugs of the other type.  (The latter are needed to fill the holes in the arms at those positions.)  Three plugs are needed because of the staggered heights of neighboring tiles.  If a prototile that is a topological sphere is desired, the plugs can be moved toward the middle of their respective faces so that the left and right side plugs meet and the holes in the arms are converted to U-shaped slots. 

To fill 3D space, the staggered layer can be stacked.  Note that in its current form the tiles in a single stacked column do not have to have identical orientations; the 2D tilings constituting successive staggered layers need not coincide, though each must be a version of the nonperiodic tiling.  If desired, a bump could be placed on the  segments on the top face of the tile directly over the point where purple and black stripes cross, with a matching indentation on the bottom face, so the tiling would be unique (periodic in the stacking direction).

Thus we see that matching rules equivalent to those of the 2D tile can be enforced by the shape of a simply connected {\em three--dimensional} prototile.
\begin{figure}[tb]
\begin{center}
\includegraphics[scale=0.21]{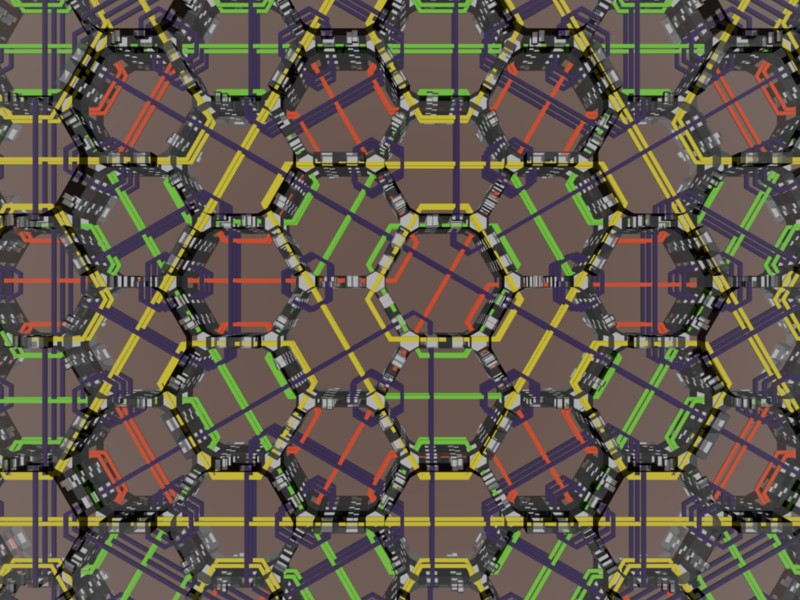}
\includegraphics[scale=0.21]{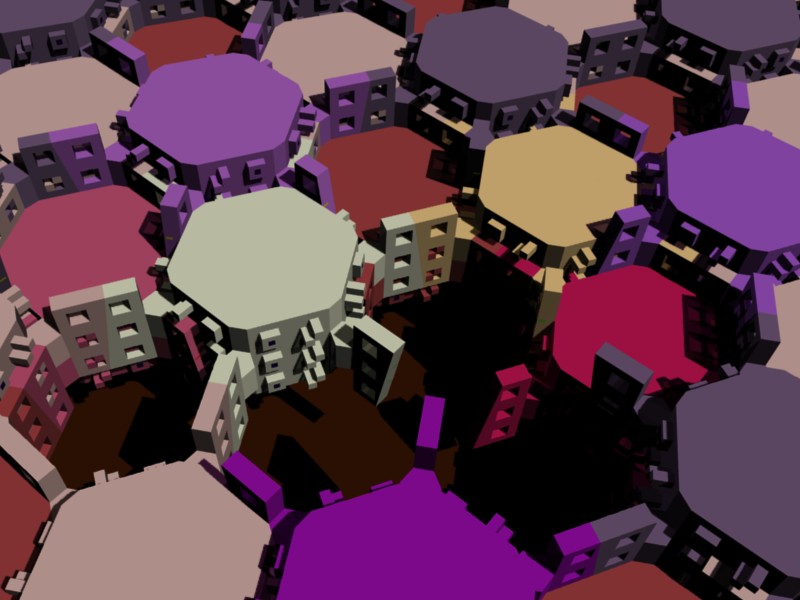}
\caption{Two renderings of portions of the 3D tiling.  For visual clarity, the purple stripes on translucent tiles at different heights are rendered in different colors.}
\label{fig3Dtiling}
\end{center}
\end{figure}
The space--filling tiling forced purely by the shape of this tile consists of a corrugated slab isomorphic to the structure forced by the 2D tiles, as shown in Figure~\ref{fig3Dtiling}, which may be stacked periodically to fill the 3D space.   This forced structure means that the tile satisfies a rather strict definition of an einstein --- the strictest definition currently known to be satisfiable.  Though the periodicity in the third dimension makes this a weakly nonperiodic tiling by Goodman-Strauss's definition,  it has a very different character from the SCD type of weak nonperiodicity.   In particular, the 3D tiling does satisfy our definition of nonperiodic, which indicates that the structure has complex correlations over large scales.

Our 3D prototile, like the SCD prototile, is not isomorphic to its reflection.  It is not possible, however, to construct a tiling that contains a mixture of the two enantiomorphs.  The shapes of the plugs enforcing the black stripe rule do not allow placing a left--handed tile adjacent to a right--handed one.  Thus the prototile is aperiodic even if one does not explicitly prohibit reflections.  For the SCD rhombic biprism, which is also chiral, reflections allow the construction of a periodic tiling with layers that alternately twist clockwise and counterclockwise.  In order to prevent mixing of the two enantiomorphs, one may decorate the SCD tile with chiral plugs, though the prototile then loses the appealing property of convexity.

\section*{The einstein pattern}
The fundamental structure of the tiling is visually evident in the patterns of black rings and purple rings in Figure~\ref{figpartialtrans}.  The black rings form truncated triangles with side lengths related by powers of two.  Any set of triangles of the same size forms a periodic honeycomb pattern equivalent to that of the smallest ones.  As expected from the similarity of the black and purple decorations of Figure~\ref{figblackpurple}, the purple rings form exactly the same pattern, rotated by $\pi/2$ and scaled up by a factor of $\sqrt{3}$.

Two proofs of aperiodicity are given in \cite{SocolarTaylorJCT}.  The first begins with an inspection of the possible ways of surrounding a given tile, which quickly reveals that a subset consisting of $3/4$ of all the tiles in the plane must be arranged to form the honeycomb of smallest black rings.  One then shows that the markings of those tiles induce precisely the same set of rules applied to the remaining tiles, so that $3/4$ of those will have to form the truncated vertices of the next largest honeycomb of black triangles.  Iterating the reasoning implies that there is no largest honeycomb, so that for any finite $r$, there will be a finite density of tiles that do not participate in partial translational symmetries with spacings smaller than $r$.  

\begin{figure}[tb]
\begin{center}
\includegraphics[scale=1.35]{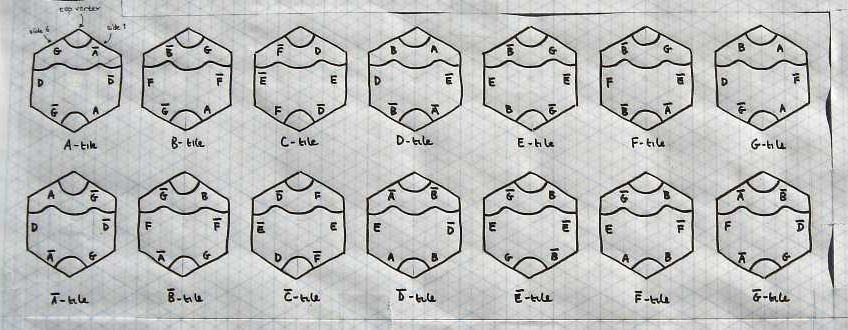}

\vspace{6pt}

\includegraphics[scale=0.32]{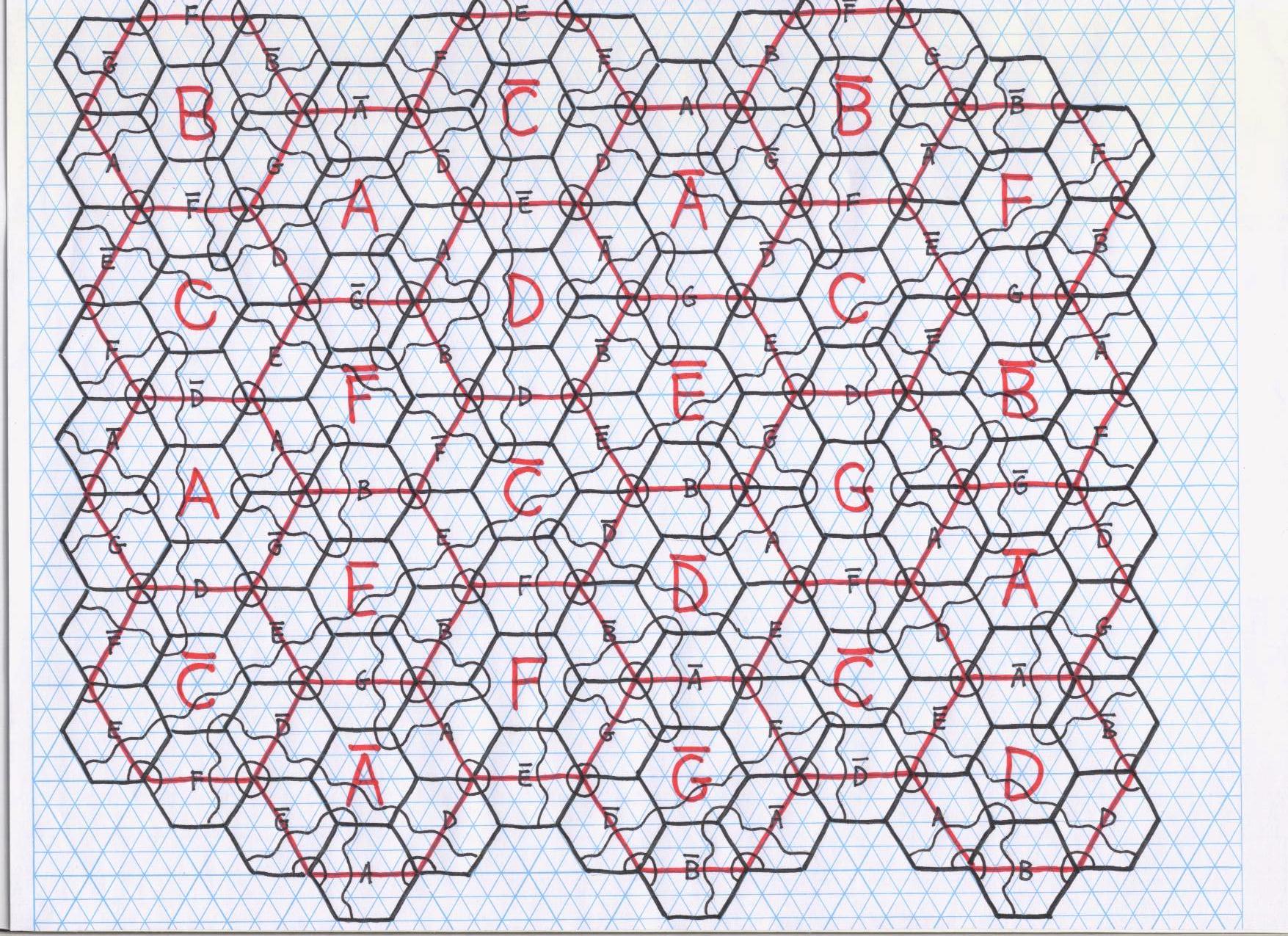}
\caption{Top: The fourteen tile set.  Bottom:  Illustration of the substitution rule.   The label of the central black hexagon within each red hexagon is $C$ if the red label is unbarred and $\bar{C}$ if the red label is barred.}
\label{figtaylorpic}
\end{center}
\end{figure}

The second makes use of the invariance of the tiling under inflation, a procedure in which tiles are grouped into larger tiles that obey the same matching rules on the larger scale.
By identifying seven distinct local environments of each chiral tile type (fourteen environments altogether) and assigning central tiles in them labels $A$ through $G$ and $\bar{A}$ through $\bar{G}$, it is possible to obtain the tiling from an iterated substitution rule, as shown in Figure~\ref{figtaylorpic}.   The scale factor associated with the substitution rule is 2, which implies that the tiling is limit--periodic (rather than quasiperiodic).~\cite{Gahler97,Godreche89}  A proof that the pattern of tile types can be enforced by a single prototile is given in~\cite{SocolarTaylorJCT}.

A curious feature of the set of forced tilings is that there is a particular arrangement of three tiles around a vertex for which the entire tiling is uniquely determined; i.e., there is a local configuration that has a unique extension to the whole plane.~\cite{SocolarTaylorJCT}  
This may be surprising, as one might expect the uniqueness of the extension to imply that the tiling must be periodic.  
Almost every finite patch that appears in a complete tiling appears an infinite number of times and permits an infinite number of distinct extensions to the entire plane.  There is, however, one particular tiling (plus its mirror image) that contains a single 3-fold symmetric vertex which does not appear in any of the other tilings.  The situation is analogous to having a decapod defect at the center of a Penrose tiling~\cite{Gardner77}, but the ``defect'' in the present case does not violate the matching rules in its interior.  

\begin{figure}[tb]
\begin{center}
\includegraphics[scale=0.85]{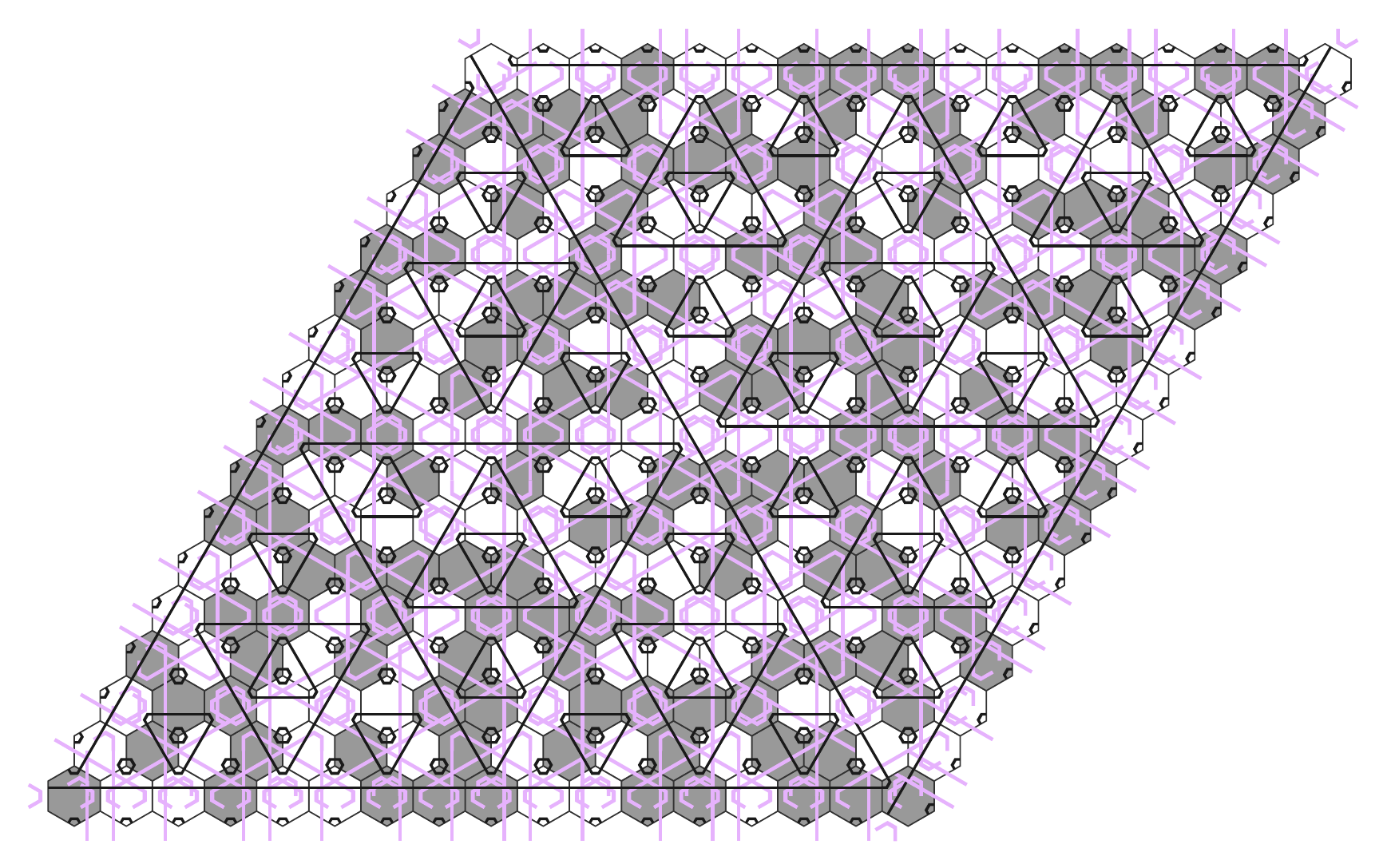}
\caption{Another coloring of the forced 2D tiling.  The purple stripes have been moved outward from the center of each tile to make it easier to see the purple triangles and the handedness of each tile is encoded in the gray and white tile colors.}
\label{figparity}
\end{center}
\end{figure}

\begin{figure}[tb]
\begin{center}
\includegraphics[scale=0.75]{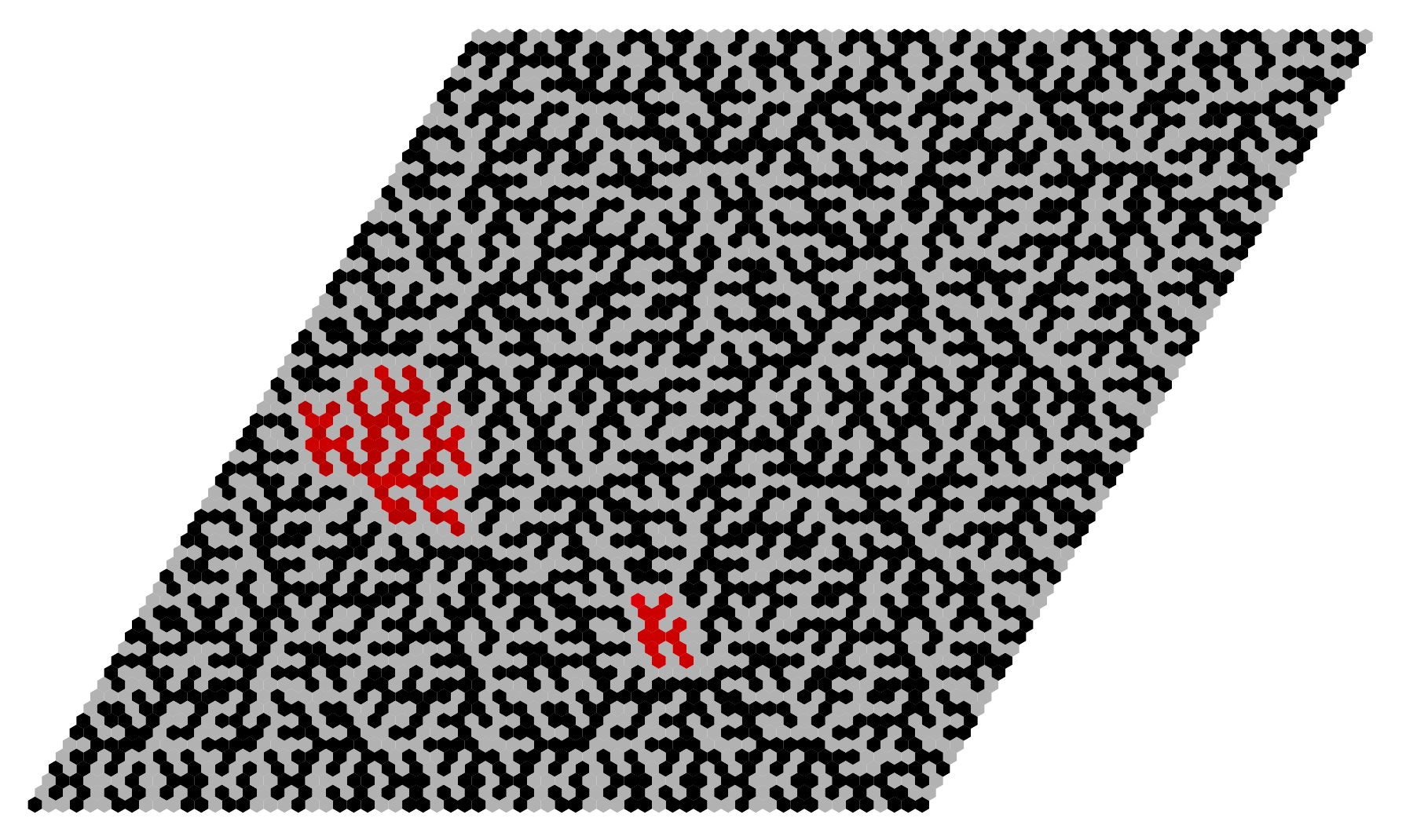}
\caption{The pattern of mirror image tiles.  The tile at the lower right corner of the rhombus corresponds to $(i,j) = (0,0)$.  Red tiles show the islands mentioned in the text.}
\label{figislands}
\end{center}
\end{figure}

A good visualization of the complexity of the tiling is obtained by shading the two mirror images differently,
as shown in Figure~\ref{figparity}. 
Figure~\ref{figislands} shows a larger portion of the tiling with one tile type shaded light grey and the mirror image tile shaded dark grey.  We have noticed a curious feature of this pattern.  There are islands of 13 dark (or light) tiles that are surrounded completely by light (or dark) tiles.  We refer to these as ``llamas.''  (See Figure~\ref{figislands}.)  Islands of 63 tiles can also be seen in the figure.  These are obtained from llamas and a few nearby tiles by application of the substitution rule.  We have also found islands of 242 tiles formed by a second iteration of the substitution rule.  Those islands each surround one llama, so the total size of the patch is 255 tiles.  We have not determined, however, whether islands of arbitrarily large size exist or whether the fraction of tiles that are {\em not} in an island of some finite size is nonzero.

The parity pattern can be specified completely as a function of tile locations with a closed--form expression.~\cite{SocolarTaylorJCT}.  An unexpected feature is the emergence of regular paperfolding sequences (A014577 of the Online Encyclopedia of Integer Sequences).  Inspection of the substitution rules for hexagons along certain rays shows that they recapitulate precisely the iterative rule that produces the regular paperfolding sequence.  This connection suggests that the full 2D pattern exhibits a rich algebraic structure that holds additional surprises and possibly affords a new window onto the properties of paperfolding and related sequences.

\section*{Closing remarks}
We have exhibited a tile that lies in a distinct new class --- a single tile that forces nonperiodicity in a space--filling tiling --- and we have presented a supporting classification scheme that captures certain intuitive distinctions between classes of nonperiodic tilings.

If mirror image tiles are counted as equivalent to the original tile, and if disconnected tiles or tiles with cutpoints are allowed, we have a 2D tile that forces a nonperiodic tiling (in exactly the same sense that the Robinson tilings are nonperiodic).  Our 3D construction gives the long sought simply connected einstein with matching rules enforced by shape alone (and no mirror image tile required).  The structure of our 3D aperiodic tile is somewhat complex and does not appear open to simplification, but two elements of the construction suggest new directions in the search for an einstein.  First, we use the possibility of rotation in 3D to create a single tile that is equivalent to two different 2D tiles, mirror images in the present case.  Second, we utilize the third dimension to implement rules that require either disconnected tiles or color matching in 2D.

The crucial point, in our view, is that our 3D prototile is the first known to force a nonperiodic structure that cannot be easily anticipated by examination of a single tile.  It is still interesting, however, to search for a single, simply connected 2D or 3D prototile that forces {\em maximal} nonperiodicity by shape alone, or one that does not permit any weakly nonperiodic tilings.  As we write, we are aware of several current computer--based searches for an aperiodic topological disk in 2D.  The general strategy is to enumerate all possible prototiles consisting of the union of simple triangles (polyiamonds), squares (polyominoes), hexagons (polyhexes), or certain pairs of triangles (polykleins), and, for each prototile, to examine all possible ways of forming small portions of a space--filling tiling.  Typically, one quickly finds a portion that can tile the whole plane periodically, or one finds that the prototile does not admit a space--filling tiling at all.  In each case, there is a measure of complexity associated with how many prototiles are needed to form the unit cell (the smallest anisohedral number permitted) or how many rings of tiles can be added around a central one before an irresolvable conflict is encountered (the Heesch number).  \cite{GoodmanStraussQs, Socolar07,isohedralweb,heeschweb}  The hope is that these computer searches will turn up a prototile that does not appear to have a finite Heesch number but does appear to have an infinite isohedral number.  Such a prototile would then have to be examined analytically to establish that it really does have both properties.  To date, the largest anisohedral number discovered is 10, which is achieved for a particular 16--hex~\cite{isohedralweb}. It is not yet clear whether computer search will beat human creativity to finding the elusive unmarked, simply connected, two-dimensional einstein --- if such a thing exists at all.

\section*{Acknowledgments}
We thank Chaim Goodman--Strauss, Marjorie Senechal, and an anonymous referee for helpful comments on earlier drafts.

\end{document}